\documentclass[a4paper]{amsart}
\usepackage{eucal,latexsym}

\newlength{\xx}
\newlength{\xd}
\DeclareMathOperator\id{id}
\DeclareMathOperator\im{im}
\DeclareMathOperator\Lin{Lin}
\DeclareMathOperator\Hom{Hom}

\DeclareMathOperator\SL{SL}
\DeclareMathOperator\SO{SO}
\DeclareMathOperator\SP{Sp}
\DeclareMathOperator\Sph{S}

\theoremstyle{definition}
\newtheorem{theorem}{Theorem}
\newtheorem*{lemma}{Lemma}

\newcommand\de{{\rm d}}
\newcommand\dd{{\rm d}'\hspace{-1pt}}
\newcommand\e{\varepsilon}
\newcommand\ox{\otimes}

\newcommand\dC{{\mathbf C}}
\newcommand\dP{{\mathbf P}}

\newcommand\cA{{\mathcal A}}
\newcommand\cC{{\mathcal C}}
\newcommand\cB{{\mathcal B}}
\newcommand\cL{{\mathcal L}}
\newcommand\cM{{\mathcal M}}
\newcommand\cO{{\mathcal O}}
\newcommand\cR{{\mathcal R}}
\newcommand\cV{{\mathcal V}}
\newcommand\cW{{\mathcal W}}

\begin{document}

\thispagestyle{empty}
\vspace*{1.6cm}

\begin{center}
DERIVATIONS WITH QUANTUM GROUP ACTION \\
\mbox{} \\
\mbox{} \\
Ulrich Hermisson \\
\mbox{} \\
Fachbereich Mathematik, Universit\"at Leipzig \\
Augustusplatz 10, 04109 Leipzig, Germany \\
uhermiss@rz.uni-leipzig.de \\
\mbox{} \\
\mbox{} \\
\mbox{} \\
\mbox{} \\
\parbox[t]{9.6cm}
{Abstract. The derivations of a left coideal subalgebra $\cB$ of a Hopf
algebra $\cA$ which are compatible with the comultiplication of $\cA$ (that
is, the covariant first order differential calculi, as defined by
Woronowicz, on a quantum homogeneous space) are related to certain right
ideals of $\cB$. The correspondence is one-to-one if $\cA$ is faithfully
flat as a right $\cB$-mod\-ule. This generalizes the result for $\cB=\cA$
due to Woronowicz. A definition for the dimension of a first order
differential calculus at a classical point is given. For the quantum
2-sphere $\Sph^2_{qc}$ of \mbox{Podle\'s} under the assumptions
$q^{n+1}\neq 1$ and $c\neq -q^{2n}/(q^{2n}+1)^2$ for all
$n=0,\,1,\,\ldots\,$, three 2-di\-men\-sional covariant first order
differential calculi exist if $c=0$, one exists if $c=\mp q/(\pm q+1)^2$ and
none else. This extends a result of Podle\'s.} \\
\mbox{} \\
\mbox{} \\
\mbox{} \\
\mbox{}
\end{center}

\section{\rm PRELIMINARIES}

A derivation of an algebra $\cB$ over $\dC$ (the complex numbers) is defined
as a $\dC$-lin\-e\-ar map $\de$ from $\cB$ into a $\cB$-bi\-mod\-ule
satisfying the Leibniz rule
\begin{displaymath}
\de(a\,b) = a\,\de b + \de a\,b \qquad \mbox{for all } a,\,b\in\cB.
\end{displaymath}
In this paper, $\de a\,b$ means $(\de a)\,b$. We set
$\Gamma(\de) = \Lin_\dC\{a\,\de b\;|\;a,\,b\in\cB\}$ (the $\dC$-lin\-e\-ar
span). We write $\dd\leq\de$, if $\dd$ and $\de$ are derivations of $\cB$
and the $\dC$-lin\-e\-ar map
$\Gamma(\de)\to\Gamma(\dd):a\,\de b\mapsto a\,\dd b$ is well-defined, and
consider derivations $\de$, $\dd$ of $\cB$ identical, if $\dd\leq\de$ and
$\de\leq\dd$. The set of derivations of $\cB$ with $\leq$ is a
complete lattice, this follows from \cite{Wo} Prop.~1.1.
If $\cB$ is a $*$-al\-ge\-bra, then $a\,\de^*b := \de(b^*)\,a^*$
defines an involution on the set of derivations of $\cB$.

We denote by $\cA$ a Hopf algebra over $\dC$ with comultiplication $\Delta$,
counit~$\e$ and antipode $S$, cf.~\cite{Sw}. We set $\otimes=\otimes_\dC$
and
\begin{displaymath}
a_{(1)}\ox a_{(2)} = \Delta(a),
\quad a_{(1)}\ox\cdots\ox a_{(n+1)}
= a_{(1)}\ox\cdots\ox a_{(n-1)}\ox\Delta(a_{(n)})
\end{displaymath}
for $n = 2,\,3,\,\ldots\,$ (Sweedler's notation) and use the map
$^+:\cA\to\cA$ defined by $a^+ = a - \e(a)\,1$. We assume that $\cB$ is a
subalgebra of $\cA$ and a left coideal,
i.e.~$\Delta(\cB)\subseteq\cA\ox\cB$, and call a derivation $\de$ of $\cB$
equivariant if and only if the $\dC$-lin\-e\-ar map
\begin{displaymath}
\Gamma(\de)\to\cA\ox\Gamma(\de):
a\,\de b\mapsto a_{(1)}\,b_{(1)}\ox a_{(2)}\,\de b_{(2)},
\qquad a,\,b\in\cB,
\end{displaymath}
is well-defined. The notion of (equivariant) derivation is the same as that
of (covariant) first order differential calculus introduced in \cite{Wo},
\cite{Po2}.

The algebra $\cB$ may be viewed as the function algebra of a quantum
homogeneous space associated to the quantum group with the function algebra
$\cA$. Accordingly, as proposed in \cite{Wo} and \cite{Po2}, \cite{Po3}, the
equivariant derivations may be considered constituents of flexibilized
(deformed) laws of nature with differential operations which are supposed to
be still invariant under the quantum group action.

In Section \ref{otocorr}, a way is prepared for determining the equivariant
derivations e.g.~for the quantizations of symmetric spaces in \cite{NS},
\cite{Di}. It is used in Section \ref{classif} in the case of the quantum
2-sphere of Podle\'s \cite{Po1} for classifying the 2-di\-men\-sional
covariant first order differential calculi. Their existence is proved by
construction in Section \ref{constr}.

\section{\rm ONE-TO-ONE CORRESPONDENCES}
\label{otocorr}

\begin{theorem}
\label{th1}
Let $\cA$ be a Hopf algebra, $\cB$ a left coideal subalgebra of
$\cA$. \\*[\xx]
(i) Let $\cR$ be a right ideal of $\cB^+$. Let $p:\cB^+\!\to\cB^+\!/\cR$ be
the canonical projection. Then $a\,\de b := a\,b_{(1)}\ox p(b_{(2)}^+)$
uniquely determines an equivariant derivation $\de_\cR:=\de$ of $\cB$. Let
$\overline{\hspace{0.75em}\rule{0pt}{1.25ex}}:\cA\to\cA/(\cB^+\!\cA)$ be the
canonical projection, $\overline{\Delta} :=
(\overline{\hspace{0.75em}\rule{0pt}{1.25ex}}\ox\id)\circ\Delta$. Then
$\cR' := \overline{\Delta}^{\,-1}(\overline{\cA}\ox\cR)$ is a right ideal of
$\cB^+$, $\cR'\subseteq\cR$, such that
$\overline{\Delta}(\cR')\subseteq\overline{\cA}\ox\cR'$ and
$\de_\cR=\de_{\cR'}$. \\[\xx]
(ii) Let $\de$ be a derivation of $\cB$. Then $\cR_\de :=
\bigl\{\sum_i\,\e(a_i)\,b_i^+\;\big|\;\sum_i\,a_i\,\de b_i = 0\bigr\}$ is a
right ideal of $\cB^+$. If $\de$ is equivariant, then
$\overline{\Delta}(\cR_\de)\subseteq\overline{\cA}\ox\cR_\de$. \\[\xx]
(iii) The maps $\cR\mapsto\de_\cR$, $\de\mapsto\cR_\de$ establish a
one-to-one correspondence of \\*[\xx]
$\/\quad\bullet$
$\{\cR_\de\;|\;\de\mbox{ an equivariant derivation of }\cB\}$ and \\*[\xx]
$\/\quad\bullet$
$\{\de_\cR\;|\;\cR\mbox{ a right ideal of }\cB^+\}$, \\[\xx]
where $\leq$ for the derivations corresponds to $\supseteq$ for the right
ideals. Furthermore, $\cR_{\de_\cR}\subseteq\cR$, if $\cR$ is a right ideal
of $\cB^+$, and $\de_{\cR_\de}\leq\de$, if $\de$ is an equivariant
derivation of $\cB$.
\end{theorem}

\begin{proof}
(i) The left module operation on
$\Gamma(\de) = \Lin_\dC\{a\,\de b\;|\;a,\,b\in\cB\}$ is determined by
$c\,(a\,\de b) = (c\,a)\,\de b$ and the right module operation by
$(a\,\de b)\,c = a\,\de(b\,c) - a\,b\,\de c$. This proves uniqueness. To
prove existence, we must show that the right module operation is
well-defined (a) and satisfies the right module axioms (b), furthermore
that $\de$ is equivariant (c). Well-definedness and axioms of the left
module operation, the bimodule axiom and the Leibniz rule clearly hold
true. \\[\xx]
(a) If $\sum_i\,a_i\,\de b_i = 0$, that is,
$\sum_i\,a_i\,b_{i(1)}\ox p(b_{i(2)}^+) = 0$, then
\begin{displaymath}
\begin{array}{l}
\bigl(\sum_i\,a_i\,\de b_i\bigr)\,c \\*[\xx]
\quad = \sum_i\,(a_i\,\de(b_i\,c) - a_i\,b_i\,\de c) \\[\xx]
\quad = \sum_i\,(a_i\,b_{i(1)}\,c_{(1)}\ox p((b_{i(2)}\,c_{(2)})^+)
- a_i\,b_i\,c_{(1)}\ox p(c_{(2)}^+)) \\[\xx]
\quad = \sum_i\,(a_i\,b_{i(1)}\,c_{(1)}\ox
p(b_{i(2)}\,c_{(2)} - \e(b_{i(2)}\,c_{(2)})) \\[\xx]
\quad\qquad - a_i\,b_{i(1)}\,c_{(1)}\ox
p(\e(b_{i(2)})\,c_{(2)} - \e(b_{i(2)}\,c_{(2)}))) \\[\xx]
\quad = \sum_i\,a_i\,b_{i(1)}\,c_{(1)}\ox p(b_{i(2)}^+\,c_{(2)}) \\*[\xx]
\quad = \bigl(\sum_i\,a_i\,b_{i(1)}\ox p(b_{i(2)}^+)\bigr)
\,(c_{(1)}\ox c_{(2)}) = 0,
\end{array}
\end{displaymath}
since $\cR$ is a right ideal of $\cB$. Hence, the specified right module
operation is well-defined. \\[\xx]
(b) We calculate that
\begin{displaymath}
\begin{array}{rl}
((a\,\de b)\,c)\,d & = (a\,\de(b\,c))\,d - (a\,b\,\de c)\,d \\*[\xx]
& = a\,\de(b\,c\,d) - a\,b\,c\,\de d - a\,b\,\de(c\,d) + a\,b\,c\,\de d
= (a\,\de b)\,(c\,d)
\end{array}
\end{displaymath}
and, since $\de 1 = 1\ox p(1^+) = 0$,
$(a\,\de b)\,1 = a\,\de(b\,1) - a\,b\,\de 1 = a\,\de b$. \\[\xx]
(c) The map $a\,\de b\mapsto a_{(1)}\,b_{(1)}\ox a_{(2)}\,\de b_{(2)}$ is
given by $\Delta\ox\id$. \\[\xx] 
We still have to show the assertions about $\cR'$. Since
$\e(\cB^+\!\cA) = \{0\}$ and
$\Delta(\cB^+\!\cA) \subseteq \cB^+\!\cA\ox\cA + \cA\ox\cB^+\!\cA$, the maps
$\e_{\overline{\cA}}:\overline{\cA}\to\dC:\overline{a}\mapsto\e(a)$ and
$\overline{\Delta}_{\overline{\cA}}:
\overline{\cA}\to\overline{\cA}\ox\overline{\cA}:
\overline{a}\mapsto\overline{a_{(1)}}\ox\overline{a_{(2)}}$ are
well-defined, and $\overline{\cA}$ with $\overline{\Delta}_{\overline{\cA}}$
is a coalgebra. If $a\in\cR'$, then
$a=\e_{\overline{\cA}}(\overline{a_{(1)}})\,a_{(2)}\in\cR$, thus
$\cR'\subseteq\cR$. Because $\cR$ is a right ideal of $\cB$,
$\overline{\Delta}(a\,b)\in\overline{\cA}\ox\cR$ for all $b\in\cB^+$, so
$\cR'$ is a right ideal of $\cB^+$. From
$\overline{a_{(1)}}\ox\overline{\Delta}(a_{(2)})
= \overline{\Delta}_{\overline{\cA}}(\overline{a_{(1)}})\ox a_{(2)}
\in\overline{\cA}\ox\overline{\cA}\ox\cR$ follows
$\overline{a_{(1)}}\ox a_{(2)}
\in\overline{\cA}\ox\overline{\Delta}^{\,-1}(\overline{\cA}\ox\cR)
= \overline{\cA}\ox\cR'$, thus
$\overline{\Delta}(\cR')\subseteq\overline{\cA}\ox\cR'$. Let
$\sum_i\,a_i\,\de b_i = 0$ with $\de=\de_\cR$, that is,
$\sum_i\,a_i\,b_{i(1)}\ox b_{i(2)}^+\in\cA\ox\cR$. Then
\begin{displaymath}
\begin{array}{l}
\sum_i\,a_i\,b_{i(1)}\ox\overline{\Delta}(b_{i(2)}^+) \\*[\xx]
\quad = \sum_i
\,\bigl(a_{i(1)}\,b_{i(1)}\ox\overline{\e(a_{i(2)})\,b_{i(2)}}\ox b_{i(3)}^+
+ a_i\,b_{i(1)}\ox\overline{b_{i(2)}^+}\ox 1\bigr) \\*[\xx]
\quad = \sum_i
\,a_{i(1)}\,b_{i(1)}\ox\overline{a_{i(2)}\,b_{i(2)}}\ox b_{i(3)}^+
\in\cA\ox\overline{\cA}\ox\cR.
\end{array}
\end{displaymath}
This implies that $\sum_i\,a_i\,b_{i(1)}\ox b_{i(2)}^+
\in\cA\ox\overline{\Delta}^{\,-1}(\overline{\cA}\ox\cR) = \cA\ox\cR'$,
therefore $\sum_i\,a_i\,\de b_i = 0$ with $\de=\de_{\cR'}$. Hence,
$\de_{\cR'}\leq\de_\cR$, while $\de_\cR\leq\de_{\cR'}$ follows from
$\cR'\subseteq\cR$. \\[\xx]
(ii) If $\sum_i\,a_i\,\de b_i = 0$, then $\sum_i\,a_i\,\de b_i\,c
= \sum_i\,(a_i\,\de(b_i\,c) - a_i\,b_i\,\de c) = 0$ and
\begin{displaymath}
\begin{array}{l}
\sum_i\,(\e(a_i)\,(b_i\,c)^+ - \e(a_i\,b_i)\,c^+) \\*[\xx]
\quad = \sum_i\,(\e(a_i)\,b_i\,c - \e(a_i)\,\e(b_i\,c)
- \e(a_i\,b_i)\,c + \e(a_i\,b_i)\,\e(c)) \\*[\xx]
\quad = \sum_i\,\e(a_i)\,b_i^+\,c \in\cR_\de.
\end{array}
\end{displaymath}
Therefore, $\cR_\de$ is a right ideal of $\cB^+$. If $\de$ is equivariant
and $\sum_i\,a_i\,\de b_i = 0$, then
$\sum_i\,f(a_{i(1)}\,b_{i(1)})\,a_{i(2)}\,\de b_{i(2)} = 0$ for any
$\dC$-lin\-e\-ar functional $f$ on $\cA$. This implies that
$\sum_i\,f(a_i\,b_{i(1)})\,b_{i(2)}^+\in\cR_\de$. From this and
\begin{displaymath}
\begin{array}{rl}
\overline{\Delta}\bigl(\sum_i\,\e(a_i)\,b_i^+\bigr)
& = \sum_i\,\bigl(\overline{\e(a_i)\,b_{i(1)}}\ox b_{i(2)}
- \overline{\e(a_i\,b_i)\,1}\ox 1\bigr) \\*[\xx]
& = \sum_i\,\bigl(\overline{a_i\,b_{i(1)}}\ox b_{i(2)}
- \overline{a_i\,b_i}\ox 1\bigr) \\*[\xx]
& = \sum_i\,\overline{a_i\,b_{i(1)}}\ox b_{i(2)}^+
\end{array}
\end{displaymath}
we conclude that
$\overline{\Delta}(\cR_\de)\subseteq\overline{\cA}\ox\cR_\de$. \\[\xx]
(iii) Directly from the definitions, we have \\*[\xx]
(a) $(\cR'\subseteq\cR)\Rightarrow(\de_\cR\leq\de_{\cR'})$ for right ideals
$\cR$, $\cR'$ of $\cB^+$, and \\*[\xx]
(b) $(\dd\leq\de)\Rightarrow(\cR_\de\subseteq\cR_{\dd})$ for derivations
$\de$, $\dd$ of $\cB$. \\[\xx]
We show that, in addition, \\*[\xx]
(c) $\cR_{\de_\cR}\subseteq\cR$, if $\cR$ is a right ideal of
$\cB^+$, and \\*[\xx]
(d) $\de_{\cR_\de}\leq\de$, if $\de$ is an equivariant derivation of
$\cB$. \\[\xx]
(c) Let $\sum_i\,a_i\,\de b_i = 0$ with $\de=\de_\cR$, that is,
$\sum_i\,a_i\,b_{i(1)}\ox p(b_{i(2)}^+) = 0$ with the canonical projection
$p:\cB^+\!\to\cB^+\!/\cR$. Application of $\e\ox\id$ leads to
$p\bigl(\sum_i\,\e(a_i)\,b_i^+\bigr) = 0$, therefore
$\sum_i\,\e(a_i)\,b_i^+\in\cR$. \\[\xx]
(d) If $\sum_i\,a_i\,\de b_i = 0$, then
$\sum_i\,f(a_{i(1)}\,b_{i(1)})\,a_{i(2)}\,\de b_{i(2)} = 0$ for any
$\dC$-lin\-e\-ar functional $f$ on $\cA$, since $\de$ is equivariant.
Thus $\sum_i\,f(a_i\,b_{i(1)})\,b_{i(2)}^+\in\cR_\de$ and
$\sum_i\,f(a_i\,b_{i(1)})\,p(b_{i(2)}^+) = 0$, where
$p:\cB^+\!\to\cB^+\!/\cR_\de$ is the canonical projection. This implies that
$\sum_i\,a_i\,\de_{\cR_\de}b_i
= \sum_i\,a_i\,b_{i(1)}\ox p(b_{i(2)}^+) = 0$. \\[\xx]
If $\de$ is an equivariant derivation of $\cB$, then
$\cR_{\de_{\cR_\de}} = \cR_\de$ by (c) and (d), (b). If $\cR$ is a right
ideal of $\cB^+$, then $\de_{\cR_{\de_\cR}} = \de_\cR$ by (d) and (c), (a).
This proves assertion (iii).
\end{proof}

In particular, the trivial derivation $\de_{\cB^+}$, the universal
derivation $\de_{\{0\}}$ and, if $\cA$ is commutative, the commutative
universal derivation $\de_{(\cB^+)^2}$ occur in the one-to-one
correspondence. The equivariant derivations of $\cB$ induced from those of
$\cA$ and e.g.~the calculi in \cite{AS}, \cite{Her}, \cite{Po2}, \cite{Po3}
also correspond to right ideals of $\cB^+$ in this way. For $\cB=\cA$,
Woronowicz \cite{Wo} shows that all equivariant derivations have this
property. We generalize this result using a theorem of Takeuchi \cite{T}
which requires further notations.

Let $\cC$ be a coalgebra, $\cW$ a right $\cC$-co\-mod\-ule and $\cV$ a left
$\cC$-co\-mod\-ule. The cotensor product $\cW\,\Box_\cC\,\cV$ is defined as
the subspace
\begin{displaymath}
\textstyle \bigl\{\sum_i\,w_i\ox v_i\in\cW\ox\cV\;\big|\;
\sum_i\,w_{i(1)}\ox w_{i(2)}\ox v_i
= \sum_i\,w_i\ox v_{i(1)}\ox v_{i(2)}\bigr\}
\end{displaymath}
of $\cW\ox\cV$; the Sweedler notation is used for the $\cC$-co\-mod\-ule
operations $\cW\to\cW\ox\cC$ and $\cV\to\cC\ox\cV$. The category of the left
$\cB$-mod\-ules $\cM$ with left $\cA$-co\-mod\-ule structure such that
$(b\,m)_{(1)}\ox(b\,m)_{(2)} = b_{(1)}\,m_{(1)}\ox b_{(2)}\,m_{(2)}$ for all
$m\in\cM$ and $b\in\cB$, together with the $\cB$-lin\-e\-ar,
$\cA$-co\-lin\-e\-ar maps, is denoted by $_\cB^\cA{\rm M}$. The category of
the left $\overline{\cA}$-co\-mod\-ules, together with the
$\overline{\cA}$-co\-lin\-e\-ar maps, is denoted by
$^{\overline{\cA}}{\rm M}$. Generally, if $\cM$ is a left $\cB$-mod\-ule,
$\overline{\hspace{0.75em}\rule{0pt}{1.25ex}}:\cM\to\cM/(\cB^+\cM)$ denotes
the canonical projection. We have shown in the proof of Theorem \ref{th1}
(i) that $\overline{\cA}$ is a coalgebra. Correspondingly, $\overline{\cM}$
is an object of $^{\overline{\cA}}{\rm M}$, if $\cM$ is an object of
$_\cB^\cA{\rm M}$. Moreover, $\cA\,\Box_{\overline{\cA}}\,\cV$ with the
induced structure of $\cA$ is an object of $_\cB^\cA{\rm M}$, if $\cV$ is an
object of $^{\overline{\cA}}{\rm M}$. That a right $\cB$-mod\-ule $\cA$ is
faithfully flat means that the functor $\cM\mapsto\cA\ox_\cB\cM$ from the
category of left $\cB$-mod\-ules to the category of $\dC$-vec\-tor spaces
preserves and reflects exact sequences.

The following result, actually the equivalent one with the opposite
multiplication and comultiplication, is contained in \cite{T}, proof of
Theorem 1. \\*[\xx]
Theorem (Takeuchi). Let $\cA$ be a Hopf algebra and $\cB$ a left coideal
subalgebra of $\cA$. If $\cA$ is faithfully flat as a right $\cB$-mod\-ule,
then the maps
\begin{displaymath}
\begin{array}{l}
\Xi:\cM\to\cA\,\Box_{\overline{\cA}}\,\overline\cM:
m\mapsto m_{(1)}\ox\overline{m_{(2)}} \qquad \mbox{and} \\*[\xx]
\Theta:\overline{\cA\,\Box_{\overline{\cA}}\,\cV}\to\cV:
\overline{\sum_i\,a_i\ox v_i}\mapsto\sum_i\,\e(a_i)\,v_i
\end{array}
\end{displaymath}
are bijective for all objects $\cM$ of $_\cB^\cA{\rm M}$ and all objects
$\cV$ of $^{\overline{\cA}}{\rm M}$. \qed

\begin{theorem}
\label{th2}
Let $\cA$ be a Hopf algebra, $\cB$ a left coideal subalgebra of
$\cA$. If $\cA$ is faithfully flat as a right $\cB$-mod\-ule, then
$\cR\mapsto\de_\cR$, $\de\mapsto\cR_\de$ as in Theorem \ref{th1}
establish a one-to-one correspondence between \\*[\xx]
$\/\quad\bullet$ the right ideals $\cR$ of $\cB^+$ with
$\overline{\Delta}(\cR)\subseteq\overline{\cA}\ox\cR$ and \\*[\xx]
$\/\quad\bullet$ the equivariant derivations of $\cB$.
\end{theorem}

\begin{proof}
If $\de$ is an equivariant derivation of $\cB$, then
$\overline{\Delta}(\cR_\de)\subseteq\overline{\cA}\ox\cR_\de$ according to
Theorem \ref{th1} (ii), so it remains to show that \\*[\xx]
(a) $\de_{\cR_\de} = \de$, if $\de$ is an equivariant derivation of
$\cB$, and \\*[\xx]
(b) $\cR_{\de_\cR} = \cR$, if $\cR$ is a right ideal of $\cB^+$ with
$\overline{\Delta}(\cR)\subseteq\overline{\cA}\ox\cR$. \\[\xx]
(a) Let $\de$ be an equivariant derivation of $\cB$. Then $\Gamma(\de)$ is
an object of~$_\cB^\cA{\rm M}$. According to Takeuchi's Theorem, the map
\begin{displaymath}
\Xi:\Gamma(\de)\to\cA\,\Box_{\overline{\cA}}\,\overline{\Gamma(\de)}:
a\,\de b\mapsto a_{(1)}\,b_{(1)}\ox\overline{a_{(2)}\,\de b_{(2)}}
= a\,b_{(1)}\ox\overline{\de b_{(2)}}
\end{displaymath}
is bijective. The kernel of
$\iota:\cB^+\!\to\overline{\Gamma(\de)}:c\mapsto\overline{\de c}$ is
$\cR_\de$: If $\sum_i\,a_i\,\de b_i=0$, then
$\overline{\de\bigl(\sum_i\,\e(a_i)\,b_i^+\bigr)}
= \overline{\sum_i\,\e(a_i)\,\de b_i}
= \overline{\sum_i\,a_i\,\de b_i} = 0$, and if $\overline{\de c} = 0$, then
$c_{ij},\,a_i,\,b_j\in\cB$ exist such that
$\de c - \sum_{ij}\,c_{ij}^+\,a_i\,\de b_j = 0$, therefore
$c^+\!\in\cR_\de$. The $\dC$-lin\-e\-ar map
$(\id\ox\iota^{-1})\circ\Xi:\Gamma(\de)\to\Gamma(\de_{\cR_\de}):
a\,\de b\mapsto a\,\de_{\cR_\de} b$ is injective and surjective, thus
$\de_{\cR_\de}=\de$. \\[\xx]
(b) Let $\cR$ be a right ideal of $\cB^+$ with
$\overline{\Delta}(\cR)\subseteq\overline{\cA}\ox\cR$. Then $\cB^+\!/\cR$
is an object of $^{\overline{\cA}}{\rm M}$. Furthermore
$\Gamma(\de_\cR)\subseteq\cA\,\Box_{\overline{\cA}}\,(\cB^+\!/\cR)$, which
follows from $a\,\de_\cR b = a\,b_{(1)}\ox p(b_{(2)}^+)$, where
$p:\cB^+\!\to\cB^+\!/\cR$ is the canonical projection, and
$a_{(1)}\,b_{(1)}\ox\overline{a_{(2)}\,b_{(2)}}\ox p(b_{(3)}^+)
= a\,b_{(1)}\ox\overline{b_{(2)}}\ox p(b_{(3)}^+)$. According to Takeuchi's
Theorem, the map
\begin{displaymath}
\Theta\big|_{\overline{\Gamma(\de_\cR)}}:
\overline{\Gamma(\de_\cR)}\to\cB^+\!/\cR:
\overline{a\,\de_\cR b} =
\overline{a\,b_{(1)}\ox p(b_{(2)}^+)}\mapsto\e(a)\,p(b^+)
\end{displaymath}
is injective. The kernel of
$\iota:\cB^+\!\to\overline{\Gamma(\de_\cR)}:c\mapsto\overline{\de_\cR c}$ is
$\cR_{\de_\cR}$, see the proof of (a), thus
$\ker(\Theta\circ\iota) = \cR_{\de_\cR}$. Since $\Theta\circ\iota = p$, we
obtain $\cR_{\de_\cR}=\cR$.
\end{proof}

In Section \ref{constr}, we give examples of equivariant derivations
which do not arise from a right ideal as in Theorem \ref{th1} (i), so
the statement of the theorem without the condition of faithful
flatness is false. However, due to M\"uller and Schneider \cite{MS}
this condition is verified for the quantizations of symmetric spaces
by Noumi, Dijkhuizen and Sugitani \cite{NS}, \cite{Di} and for the
quantized flag manifolds \cite{SD}.

\section{\rm CLASSIFICATION}
\label{classif}

We call $\dim_{\e,{\rm l}}\de := \dim_\dC(\Gamma(\de)/(\cB^+\Gamma(\de)))$
the left dimension and analogously
$\dim_{\e,{\rm r}}\de := \dim_\dC(\Gamma(\de)/(\Gamma(\de)\cB^+))$ the
right dimension of a first order differential calculus $\de$ over $\cB$ at
the classical point $\e$. If $\cB$ is the algebra of regular functions on a
nonsingular affine algebraic variety, then
$\dim_{\e,{\rm l}}\de_{(\cB^+)^2} = \dim_{\e,{\rm r}}\de_{(\cB^+)^2}$ is the
dimension of it. If a basis of $\Gamma(\de)$ as a left $\cB$-mod\-ule
exists, e.g.~if $\cB=\cA$ and $\de$ is left-covariant, cf.~\cite{Wo}, then
$\dim_{\e,{\rm l}}\de$ is the number of its elements.

We assume that $q\in\dC\setminus\{0\}$ and $q^n \neq 1$ for all
$n=1,\,2,\,\ldots\,$. For the quantum 2-sphere of Podle\'s we may
equivalently choose $\cA$ to be one of the quantum group function algebras
$\cO(\SL_q(2))$, $\cO(\SO_{q^2}(3))$ and $\cO(\SP_{q^{1/2}}(2))$ which are
described in \cite{RTF}. The function algebras $\cO(\Sph^2_{qc})$ of the
quantum 2-sphere, parameterized by $c\in\dC\dP^1$, are the
$\cA$-co\-mod\-ule algebras (except for one) which are isomorphic as a
comodule to the classical case and generated as an algebra by the
spin 1 subcomodule, cf.~\cite{Po1}. They are isomorphic to right
coideal subalgebras $\cB_c$ of $\cA$, so the equivalents of our
theorems with the opposite multiplication and comultiplication are
applicable (we silently assume this exchange of left and right). The
algebras $\cB_c$ are generated by three elements $e_{-1}$, $e_0$,
$e_1$ with the relations
\begin{displaymath}
\begin{array}{l}
(q^2+1)\,e_{-1}\,e_1 + e_0^2 + (q^{-2}+1)\,e_1\,e_{-1} = \rho\,1, \\*[\xx]
-q^2\,e_{-1}\,e_0 + e_0\,e_{-1} = \lambda\,e_{-1}, \\*[\xx]
(q^2+1)\,e_{-1}\,e_1 - (q^2-1)\,e_0^2 - (q^2+1)\,e_1\,e_{-1}
= \lambda\,e_0, \\*[\xx]
-q^2\,e_0\,e_1 + e_1\,e_0 = \lambda\,e_1,
\end{array}
\quad\rho,\,\lambda\in\dC,\quad(*)
\end{displaymath}
such that $c = \e(e_{-1})\,\e(e_1) : \e(e_0)^2$ and
$\Delta(e_i) = \sum_j\,e_j\ox\pi^j_i$, cf.~\cite{Po1} ($q$ is $\mu$,
$\pi^i_j$ is $d_{1,ij}$). Special values of $c$ are
$c(n) = -q^{2n}/(q^{2n}+1)^2$.

We classify the equivariant derivations $\de$ of $\cB_c$ with
$\dim_{\e,{\rm r}}\de = 2$ which arise from a left ideal as in Theorem
\ref{th1} (i). If $c\neq c(n)$ for all $n=0,\,1,\,\ldots\,$, then
$\cA$ is faithfully flat as a left $\cB_c$-mod\-ule, cf.~\cite{Mue1},
\cite{Mue2}, \cite{MS}, and according to Theorem \ref{th2} our
classification includes all equivariant derivations with
$\dim_{\e,{\rm r}}\de = 2$. We denote by
$\overline{\hspace{0.75em}\rule{0pt}{1.25ex}}$, $\overline{\Delta}$
and $\cL_\de$ the equivalents of the previously used structures with
left and right reversed. The map
$\chi:\cB_c\to\overline{\Gamma(\de)}:b\mapsto\overline{\de b}$ induces a
$\dC$-lin\-e\-ar bijection between $\cB_c^+\!/\cL_\de$ and
$\overline{\Gamma(\de)}$, see the proof of Theorem \ref{th2}, i.e.~we must
determine the left ideals $\cL$ of $\cB_c^+$ with
$\dim_\dC(\cB_c^+\!/\cL) = 2$ and
$\overline{\Delta}(\cL)\subseteq\cL\ox\overline{\cA}$. The Hochschild
coboundary maps of the quotient $\cB_c$-bi\-mod\-ule
$\overline{\Gamma(\de)}$ are defined as
\begin{displaymath}
\begin{array}{l}
\delta^0:\overline{\Gamma(\de)}\to\Hom_\dC(\cB_c,\,\overline{\Gamma(\de)}):
\overline{\omega} \mapsto
(b \mapsto b\,\overline{\omega} - \overline{\omega}\,b)
\qquad \mbox{and} \\[\xx]
\delta^n:\Hom_\dC(\cB_c^{\ox n},\,\overline{\Gamma(\de)})\to
\Hom_\dC(\cB_c^{\ox(n+1)},\,\overline{\Gamma(\de)}),
\end{array}
\end{displaymath}
$\delta^n\! f(b_0\ox\cdots\ox b_n) = b_0\,f(b_1\ox\cdots\ox b_n) \\*[\xx]
+ \sum\limits_{i=1}^n
\,(-1)^i\,f(b_0\ox\cdots\ox b_{i-1}\,b_i\ox\cdots\ox b_n)
+ (-1)^{n+1}\,f(b_0\ox\cdots\ox b_{n-1})\,b_n$ \\[\xx]
for $n=1,\,2,\,\ldots\,$. Hence, $\chi$ is a 1-co\-cy\-cle:
$a\,\chi(b) - \chi(a\,b) + \chi(a)\,b = 0$ for all $a,\,b\in\cB_c$. Let
$\tau:\cB_c\to\Hom_\dC(\overline{\Gamma(\de)},\,\overline{\Gamma(\de)})$ be
the representation of $\cB_c$ on the quotient left $\cB_c$-mod\-ule
$\overline{\Gamma(\de)}$. Using coordinates, this says
\begin{displaymath}
\begin{array}{l}
\chi_i(a\,b) =
\sum_k\,\tau_{ik}(a)\,\chi_k(b) + \chi_i(a)\,\e(b),
\qquad \chi_i(1) = 0, \\*[\xx]
\tau_{ij}(a\,b) = \sum_k\,\tau_{ik}(a)\,\tau_{kj}(b),
\qquad \tau_{ij}(1) = \delta_{ij}
\end{array}
\end{displaymath}
for all $a,\,b\in\cB_c$ and $i,\,j=1,\,2$. Given a representation $\tau$,
the 1-co\-bound\-a\-ries $\chi_i = \sum_k\,\beta_k\,\tau_{ik} - \beta_i\,\e$
with $\beta_1,\,\beta_2\in\dC$ are solutions to these equations, and further
solutions exist exactly if the first cohomology group
$(\ker\delta^1)/(\im\delta^0)$ is not $\{0\}$. The equations imply that the
functions $\chi_i$ and $\tau_{ij}$ are uniquely determined by their values
on $e_{-1}$, $e_0$, $e_1$ and exist for given values if they are compatible
with the relations $(*)$. The solutions in suitable coordinates are
$\chi_i = \sum_k\,\beta_k\,\tau_{ik} - \beta_i\,\e
+ \sum_n\,\beta'_n\,\xi^n_i$ with $\beta_i,\,\beta'_n\in\dC$, \\[\xx]
(a) $\tau(e_{-1},\,e_0,\,e_1) =
x^{-1}\,\bigl(\begin{smallmatrix}
q^{-2}\,\alpha_{-1}\,\alpha_1 & \alpha_0 \\ 0 & \alpha_{-1}\,\alpha_1
\end{smallmatrix}\bigr),
\,\bigl(\begin{smallmatrix}
\alpha_0 & -(q^2+1) \\ 0 & \alpha_0
\end{smallmatrix}\bigr),
\,x\,\bigl(\begin{smallmatrix}
q^2 & 0 \\ 0 & 1
\end{smallmatrix}\bigr)$, \\*[\xd]
$\xi^1(e_{-1},\,e_0,\,e_1) =
\bigl(\begin{smallmatrix} -\alpha_{-1} \\ 0 \end{smallmatrix}\bigr),
\,\bigl(\begin{smallmatrix} 0 \\ 0 \end{smallmatrix}\bigr),
\,\bigl(\begin{smallmatrix} \alpha_1 \\ 0 \end{smallmatrix}\bigr)$
if $x=q^{-2}\,\alpha_1$ and \\*[\xd]
$\xi^1(e_{-1},\,e_0,\,e_1) =
\bigl(\begin{smallmatrix} -1 \\ \alpha_0 \end{smallmatrix}\bigr),
\,\bigl(\begin{smallmatrix} 0 \\ -(q^2+1)\,\alpha_1 \end{smallmatrix}\bigr),
\,\bigl(\begin{smallmatrix} 0 \\ 0 \end{smallmatrix}\bigr)$
if $x=q^2\,\alpha_1$, \\[\xx]
(b) $\tau(e_{-1},\,e_0,\,e_1) =
x^{-1}\,\alpha_{-1}\,\alpha_1\,\bigl(\begin{smallmatrix}
1 & -1 \\ 0 & 1
\end{smallmatrix}\bigr),
\,\bigl(\begin{smallmatrix}
\alpha_0 & 0 \\ 0 & \alpha_0
\end{smallmatrix}\bigr),
\,x\,\bigl(\begin{smallmatrix}
1 & 1 \\ 0 & 1
\end{smallmatrix}\bigr)$, \\*[\xd]
$\xi^1(e_{-1},\,e_0,\,e_1) =
\bigl(\begin{smallmatrix}
\alpha_{-1} \\ -\alpha_{-1}
\end{smallmatrix}\bigr),
\,\bigl(\begin{smallmatrix} 0 \\ 0 \end{smallmatrix}\bigr),
\,\bigl(\begin{smallmatrix} 0 \\ \alpha_1 \end{smallmatrix}\bigr)$
if $x=\alpha_1$ and \\*[\xd]
$\xi^1(e_{-1},\,e_0,\,e_1) =
\bigl(\begin{smallmatrix} \alpha_0 \\ 0 \end{smallmatrix}\bigr),
\,\bigl(\begin{smallmatrix} -(q^2+1)\,\alpha_1 \\ 0 \end{smallmatrix}\bigr),
\,\bigl(\begin{smallmatrix} 0 \\ 0 \end{smallmatrix}\bigr)$
if $x=q^2\,\alpha_1$, \\[\xx]
(c) $\tau(e_{-1},\,e_0,\,e_1) =
\alpha_{-1}\,\alpha_1\,\bigl(\begin{smallmatrix}
x^{-1} & 0 \\ 0 & y^{-1}
\end{smallmatrix}\bigr),
\,\bigl(\begin{smallmatrix}
\alpha_0 & 0 \\ 0 & \alpha_0
\end{smallmatrix}\bigr),
\,\bigl(\begin{smallmatrix}
x & 0 \\ 0 & y
\end{smallmatrix}\bigr)$, \\*[\xd]
$\xi^1(e_{-1},\,e_0,\,e_1) =
\bigl(\begin{smallmatrix} \alpha_{-1} \\ 0 \end{smallmatrix}\bigr),
\,\bigl(\begin{smallmatrix} 0 \\ 0 \end{smallmatrix}\bigr),
\,\bigl(\begin{smallmatrix} -\alpha_1 \\ 0 \end{smallmatrix}\bigr)$
if $x=\alpha_1$ and \\*[\xd]
$\xi^1(e_{-1},\,e_0,\,e_1) =
\bigl(\begin{smallmatrix} \alpha_0 \\ 0 \end{smallmatrix}\bigr),
\,\bigl(\begin{smallmatrix} -(q^2+1)\,\alpha_1 \\ 0 \end{smallmatrix}\bigr),
\,\bigl(\begin{smallmatrix} 0 \\ 0 \end{smallmatrix}\bigr)$
if $x=q^2\,\alpha_1$, \\*[\xd]
$\xi^2(e_{-1},\,e_0,\,e_1) =
\bigl(\begin{smallmatrix} 0 \\ \alpha_{-1} \end{smallmatrix}\bigr),
\,\bigl(\begin{smallmatrix} 0 \\ 0 \end{smallmatrix}\bigr),
\,\bigl(\begin{smallmatrix} 0 \\ -\alpha_1 \end{smallmatrix}\bigr)$
if $y=\alpha_1$ and \\*[\xd]
$\xi^2(e_{-1},\,e_0,\,e_1) =
\bigl(\begin{smallmatrix} 0 \\ \alpha_0 \end{smallmatrix}\bigr),
\,\bigl(\begin{smallmatrix} 0 \\ -(q^2+1)\,\alpha_1 \end{smallmatrix}\bigr),
\,\bigl(\begin{smallmatrix} 0 \\ 0 \end{smallmatrix}\bigr)$
if $y=q^2\,\alpha_1$, \\[\xx]
where $x,\,y\in\dC\setminus\{0\}$, $\alpha_i=\e(e_i)$ and, if not specified
otherwise, $\xi^j=0$; in addition, if $c=c(1)$, \\[\xx]
(d) $\tau(e_{-1},\,e_0,\,e_1) =
\alpha_{-1}\,\alpha_1\,\bigl(\begin{smallmatrix}
x^{-1} & 0 \\ 0 & 0
\end{smallmatrix}\bigr),
\,\bigl(\begin{smallmatrix} \alpha_0 & 0 \\ 0 & 0 \end{smallmatrix}\bigr),
\,\bigl(\begin{smallmatrix}
x & 0 \\ 0 & 0
\end{smallmatrix}\bigr)$, \\*[\xd]
$\xi^1(e_{-1},\,e_0,\,e_1) =
\bigl(\begin{smallmatrix} \alpha_{-1} \\ 0 \end{smallmatrix}\bigr),
\,\bigl(\begin{smallmatrix} 0 \\ 0 \end{smallmatrix}\bigr),
\,\bigl(\begin{smallmatrix} -\alpha_1 \\ 0 \end{smallmatrix}\bigr)$
if $x=\alpha_1$ and \\*[\xd]
$\xi^1(e_{-1},\,e_0,\,e_1) =
\bigl(\begin{smallmatrix} \alpha_0 \\ 0 \end{smallmatrix}\bigr),
\,\bigl(\begin{smallmatrix} -(q^2+1)\,\alpha_1 \\ 0 \end{smallmatrix}\bigr),
\,\bigl(\begin{smallmatrix} 0 \\ 0 \end{smallmatrix}\bigr)$
if $x=q^2\,\alpha_1$, \\[\xx]
(e) $\tau(e_{-1},\,e_0,\,e_1) =
\bigl(\begin{smallmatrix} 0 & 0 \\ 0 & 0 \end{smallmatrix}\bigr),
\,\bigl(\begin{smallmatrix} 0 & 0 \\ 0 & 0 \end{smallmatrix}\bigr),
\,\bigl(\begin{smallmatrix} 0 & 0 \\ 0 & 0 \end{smallmatrix}\bigr)$ \\[\xx]
and if $c=c(2)$, \\[\xx]
(f) $\tau(e_{-1},\,e_0,\,e_1) =
\frac{q^2-1}{q^4+1}\,\alpha_0\,\bigl(\begin{smallmatrix}
0 & 1 \\ 0 & 0
\end{smallmatrix}\bigr),
\,\frac{q^2-1}{q^4+1}\,\alpha_0\,\bigl(\begin{smallmatrix}
-1 & 0 \\ 0 & q^2
\end{smallmatrix}\bigr),
\,\frac{q^2-1}{q^4+1}\,\alpha_0\,\bigl(\begin{smallmatrix}
0 & 0 \\ q^2 & 0
\end{smallmatrix}\bigr)$ \\[\xx]
and if $c=0$, \\[\xx]
(a') $\tau(e_{-1},\,e_0,\,e_1) =
x\,\bigl(\begin{smallmatrix}
q^{-2} & 0 \\ 0 & 1
\end{smallmatrix}\bigr),
\,\bigl(\begin{smallmatrix}
\alpha_0 & -(q^{-2}+1) \\ 0 & \alpha_0
\end{smallmatrix}\bigr),
\,x^{-1}\,\bigl(\begin{smallmatrix}
0 & \alpha_0 \\ 0 & 0
\end{smallmatrix}\bigr)$, \\*[\xd]
$\xi^1(e_{-1},\,e_0,\,e_1) =
\bigl(\begin{smallmatrix} 1 \\ 0 \end{smallmatrix}\bigr),
\,\bigl(\begin{smallmatrix} 0 \\ 0 \end{smallmatrix}\bigr),
\,\bigl(\begin{smallmatrix} 0 \\ 0 \end{smallmatrix}\bigr)$
if $x=q^2\,\alpha_{-1}$ and \\*[\xd]
$\xi^1(e_{-1},\,e_0,\,e_1) =
\bigl(\begin{smallmatrix} 0 \\ 0 \end{smallmatrix}\bigr),
\,\bigl(\begin{smallmatrix}
0 \\ -(q^{-2}+1)\,\alpha_{-1}
\end{smallmatrix}\bigr),
\,\bigl(\begin{smallmatrix} -1 \\ \alpha_0 \end{smallmatrix}\bigr)$
if $x=q^{-2}\,\alpha_{-1}$, \\[\xx]
(b') $\tau(e_{-1},\,e_0,\,e_1) =
\bigl(\begin{smallmatrix} 0 & s \\ 0 & 0 \end{smallmatrix}\bigr),
\,\bigl(\begin{smallmatrix}
\alpha_0 & 0 \\ 0 & \alpha_0
\end{smallmatrix}\bigr),
\,\bigl(\begin{smallmatrix}
t & 1 \\ 0 & t
\end{smallmatrix}\bigr)$ with $s\,t=0$, \\*[\xd]
$\xi^1(e_{-1},\,e_0,\,e_1) =
\bigl(\begin{smallmatrix} 0 \\ q^2\,s\,\alpha_0 \end{smallmatrix}\bigr),
\,\bigl(\begin{smallmatrix} -(q^2+1)\,s \\ 0 \end{smallmatrix}\bigr),
\,\bigl(\begin{smallmatrix} 0 \\ \alpha_0 \end{smallmatrix}\bigr)$,
if $\alpha_{-1}=0$ and $t=\alpha_1$, \\*[\xd]
$\xi^2(e_{-1},\,e_0,\,e_1) =
\bigl(\begin{smallmatrix} \alpha_0 \\ 0 \end{smallmatrix}\bigr),
\,\bigl(\begin{smallmatrix} -(q^2+1)\,\alpha_1 \\ 0 \end{smallmatrix}\bigr),
\,\bigl(\begin{smallmatrix} 0 \\ 0 \end{smallmatrix}\bigr)$,
if $\alpha_{-1}=0$ and $t=q^2\,\alpha_1$, \\[\xx]
(b'') $\tau(e_{-1},\,e_0,\,e_1) =
\bigl(\begin{smallmatrix} s & 1 \\ 0 & s \end{smallmatrix}\bigr),
\,\bigl(\begin{smallmatrix}
\alpha_0 & 0 \\ 0 & \alpha_0
\end{smallmatrix}\bigr),
\,\bigl(\begin{smallmatrix}
0 & 0 \\ 0 & 0
\end{smallmatrix}\bigr)$, \\*[\xd]
$\xi^1(e_{-1},\,e_0,\,e_1) =
\bigl(\begin{smallmatrix} 0 \\ 1 \end{smallmatrix}\bigr),
\,\bigl(\begin{smallmatrix} 0 \\ 0 \end{smallmatrix}\bigr),
\,\bigl(\begin{smallmatrix} 0 \\ 0 \end{smallmatrix}\bigr)$,
if $\alpha_1=0$ and $s=\alpha_{-1}$, \\*[\xd]
$\xi^2(e_{-1},\,e_0,\,e_1) =
\bigl(\begin{smallmatrix} 0 \\ 0 \end{smallmatrix}\bigr),
\,\bigl(\begin{smallmatrix}
-(q^{-2}+1)\,\alpha_{-1} \\ 0
\end{smallmatrix}\bigr),
\,\bigl(\begin{smallmatrix} \alpha_0 \\ 0 \end{smallmatrix}\bigr)$,
if $\alpha_1=0$ and $s=q^{-2}\,\alpha_{-1}$, \\[\xx]
(c') $\tau(e_{-1},\,e_0,\,e_1) =
\bigl(\begin{smallmatrix} s & 0 \\ 0 & u \end{smallmatrix}\bigr),
\,\bigl(\begin{smallmatrix}
\alpha_0 & 0 \\ 0 & \alpha_0
\end{smallmatrix}\bigr),
\,\bigl(\begin{smallmatrix} t & 0 \\ 0 & v \end{smallmatrix}\bigr)$
with $s\,t=u\,v=0$, \\*[\xd]
$\xi^1(e_{-1},\,e_0,\,e_1) =
\bigl(\begin{smallmatrix} 0 \\ 0 \end{smallmatrix}\bigr),
\,\bigl(\begin{smallmatrix}
-(q^{-2}+1)\,\alpha_{-1} \\ 0
\end{smallmatrix}\!\bigr),
\,\bigl(\begin{smallmatrix} \alpha_0 \\ 0 \end{smallmatrix}\bigr)$,
if $s=q^{-2}\,\alpha_{-1}$ and $t=\alpha_1$, \\*[\xd]
$\xi^2(e_{-1},\,e_0,\,e_1) =
\bigl(\begin{smallmatrix} \alpha_0 \\ 0 \end{smallmatrix}\bigr),
\,\bigl(\begin{smallmatrix}
-(q^2+1)\,\alpha_1 \\ 0
\end{smallmatrix}\bigr),
\,\bigl(\begin{smallmatrix} 0 \\ 0 \end{smallmatrix}\bigr)$,
if $s=\alpha_{-1}$ and $t=q^2\,\alpha_1$, \\*[\xd]
$\xi^3(e_{-1},\,e_0,\,e_1) =
\bigl(\begin{smallmatrix} 0 \\ 0 \end{smallmatrix}\bigr),
\,\bigl(\begin{smallmatrix}
0 \\ -(q^{-2}+1)\,\alpha_{-1}
\end{smallmatrix}\bigr),
\,\bigl(\begin{smallmatrix} 0 \\ \alpha_0 \end{smallmatrix}\bigr)$,
if $u=q^{-2}\,\alpha_{-1}$ and $v=\alpha_1$, \\*[\xd]
$\xi^4(e_{-1},\,e_0,\,e_1) =
\bigl(\begin{smallmatrix} 0 \\ \alpha_0 \end{smallmatrix}\bigr),
\,\bigl(\begin{smallmatrix}
0 \\ -(q^2+1)\,\alpha_1
\end{smallmatrix}\bigr),
\,\bigl(\begin{smallmatrix} 0 \\ 0 \end{smallmatrix}\bigr)$,
if $u=\alpha_{-1}$ and $v=q^2\,\alpha_1$, \\[\xx]
where $s,\,t,\,u,\,v\in\dC$. The solutions (b'), (c') contain (b), (c).

The condition $\overline{\Delta}(\cL)\subseteq\cL\ox\overline{\cA}$ implies
that $\chi(b_{(1)}) f(b_{(2)}) = 0$ for all $b\in\cL$ and $\dC$-lin\-e\-ar
functionals $f$ on $\cA$ with $f(\cA\cB_c^+) = \{0\}$, therefore
$\chi_1 f,\,\chi_2 f\in\Lin_\dC\{\chi_1,\,\chi_2\}$. Such a functional
(taken from \cite{DK}) is
\begin{displaymath}
f = \alpha_0\,(l^2-\e) - \alpha_1\,l^{(+)}\,l + \alpha_{-1}\,l^{(-)}\,l
\end{displaymath}
with $l=l^{(+)}_{11}$, $l^{(+)}=l^{(+)}_{12}$ and $l^{(-)}=l^{(-)}_{21}$ as
specified for $\cA=\cO(\SL_q(2))$ in \cite{RTF}. The conditions
$\chi_1 f,\,\chi_2 f\in\Lin_\dC\{\chi_1,\,\chi_2\}$, checked by evaluation
of $\chi_1$, $\chi_2$, $\chi_1 f$ and $\chi_2 f$ on $e_i$, $e_i\,e_j$ and
$e_i\,e_j\,e_k$ for $i,\,j,\,k=-1,\,0,\,1$, and
$\dim_\dC\Lin_\dC\{\chi_1,\,\chi_2\} = 2$ reduce the total number of
solutions to seven: \\[\xx]
1.+2. $c=\mp q/(\pm q+1)^2$, (a) with $x=\pm q^{-1}\,\alpha_1$ and
$\frac{\beta_1}{\beta_2}
= \frac{\alpha_1^2}{\alpha_0\,x\,(x-\alpha_1)}$. \\[\xx]
3.+4. $c=c(2)$, (f) with
$\frac{\beta_1}{\beta_2}\in\bigl\{\frac{(q^4+1)\,\alpha_1}{q^4\,\alpha_0},
\,\frac{(q^4+1)\,\alpha_1}{\alpha_0}\bigr\}$. \\[\xx]
5.--7. $c=0$, $\alpha_{-1}=0=\alpha_1$, \\*[\xd]
(b'') with $s=0$ and $\beta'_1\neq 0=\beta'_2$, \\*[\xd]
(b') with $s=t=0$ and $\beta'_1\neq 0=\beta'_2$, \\*[\xd]
(c') with $s=t=u=v=0$ and
$\beta'_1\,\beta'_4\neq\beta'_2\,\beta'_3$. \\[\xx]
5.--7. $c=0$, $\alpha_{-1}\neq 0=\alpha_1$, \\*[\xd]
(c') with $s=q^2\,\alpha_{-1}$, $u=q^4\,\alpha_{-1}$, $t=v=0$ and
$\beta_1\neq 0\neq\beta_2$, \\*[\xd]
(a') with $x=q^{-2}\,\alpha_{-1}$ and
$\frac{\beta_1}{\beta_2} = \frac{q^2}{(q^4-1)\,\alpha_0}$,
$\frac{\beta_2}{\beta'_1} = \frac{\alpha_{-1}}{(q^2-1)\,\alpha_0}$, \\*[\xd]
(c') with $s=q^{-2}\,\alpha_{-1}$, $u=q^2\,\alpha_{-1}$, $t=v=0$ and
$\frac{\beta_1}{\beta'_1} = \frac{\alpha_{-1}}{(q^2-1)\,\alpha_0}$,
$\beta_2\neq 0$. \\[\xx]
5.--7. $c=0$, $\alpha_{-1}=0\neq\alpha_1$, \\*[\xd]
(a) with $x=q^2\,\alpha_1$ and
$\frac{\beta_1}{\beta_2} = - \frac{q^2}{(q^4-1)\,\alpha_0}$,
$\frac{\beta_2}{\beta'_1} = \frac{\alpha_1}{(q^{-2}-1)\,\alpha_0}$, \\*[\xd]
(c') with $s=u=0$, $t=q^{-2}\,\alpha_1$, $v=q^{-4}\,\alpha_1$ and
$\beta_1\neq 0\neq\beta_2$, \\*[\xd]
(c') with $s=u=0$, $t=q^2\,\alpha_1$, $v=q^{-2}\,\alpha_1$ and
$\frac{\beta_1}{\beta'_2} = \frac{\alpha_1}{(q^{-2}-1)\,\alpha_0}$,
$\beta_2\neq 0$. \\[\xx]
We have carried out the calculations for all embeddings of
$\cO(\Sph^2_{qc})$ in $\cA$. For each embedding, $\e$ determines a classical
point of $\cO(\Sph^2_{qc})$, i.e.~an algebra homomorphism
$\cO(\Sph^2_{qc})\to\dC$; if $c\neq c(1)$, this is a one-to-one
correspondence. The equivariant derivations of $\cB_c$ corresponding to the
solutions 1--7 are given by $\de b\,a = \chi(b_{(1)})\ox b_{(2)}\,a$. They
are independent of the embedding and satisfy $\dim_{\e,{\rm r}}\de = 2$ for
each classical point $\e$, in the case of the solutions 1, 2 and 5--7 also
$\dim_{\e,{\rm l}}\de = 2$; all this is proved in Section \ref{constr}. The
solutions for $c=0$, not being coboundaries if $\alpha_{-1}=0=\alpha_1$, do
not correspond to derivations with an $\omega\in\Gamma(\de)$ for which
$\de a = \omega\,a - a\,\omega$ for all $a\in\cB_c$. If
$\de a = \sum_i\,\de e_i\,a_i$, then $\chi(a) = \chi(\sum_i\,\e(a_i)\,e_i)$.
Since $\dim_\dC\chi(\Lin_\dC\{e_{-1},\,e_0,\,e_1\}) = 1$ for the solutions
3, 5 and 6, but $\dim_\dC\chi(\cB_c) = 2$, these do not correspond to
derivations for which $\{\de e_{-1},\,\de e_0,\,\de e_1\}$ generates
$\Gamma(\de)$ as a right $\cB_c$-module.

\section{\rm CONSTRUCTIONS}
\label{constr}

We retain the notations of Section \ref{classif}. In the case of a quantum
group, $\cA$ is coquasitriangular, i.e.~equipped with a $\dC$-lin\-e\-ar map
$r:\cA\ox\cA\to\dC$ which satisfies
\begin{displaymath}
\begin{array}{l}
r(a_{(1)}\ox b_{(1)})\,a_{(2)}\,b_{(2)}
= b_{(1)}\,a_{(1)}\,r(a_{(2)}\ox b_{(2)}), \\*[\xx]
r(a\,b\ox c) = r(a\ox c_{(1)})\,r(b\ox c_{(2)}),
\qquad r(1\ox c) = \e(c), \\*[\xx]
r(a\ox b\,c) = r(a_{(1)}\ox c)\,r(a_{(2)}\ox b),
\qquad r(a\ox 1) = \e(a)
\end{array}
\end{displaymath}
for all $a,\,b,\,c\in\cA$. We use a construction method introduced in
\cite{Her}.

\begin{lemma}
Let $\cA$ be a coquasitriangular Hopf algebra and $\cB$ a right
$\cA$-co\-mod\-ule algebra. Let $b_1,\,\ldots,\,b_N\in\cB$ be
$\dC$-lin\-e\-arly independent elements,
$\Delta(b_i) = \sum_j\,b_j\ox\psi^j_i$ and $\nu$ a comodule algebra
endomorphism of $\cB$. Let $\Gamma$ be the $\cB$-bi\-mod\-ule generated by
the symbols $\gamma^1,\,\ldots,\,\gamma^N$ with the relations
$a\,\gamma^j = \sum_i\,\gamma^i\,\nu(a_{(1)})\,r(\psi^j_i\ox a_{(2)})$,
$a\in\cB$, $j=1,\,\ldots,\,N$. Then $\gamma^1,\,\ldots,\,\gamma^N$ is a
basis of $\Gamma$ as a right $\cB$-mod\-ule, and
$\de a := \omega\,a - a\,\omega$ with $\omega = \sum_i\,\gamma^i\,b_i$
defines an equivariant derivation of $\cB$. Moreover, if $\nu$ is bijective,
then $\gamma^i\,a
= \sum_j\,\nu^{-1}(a_{(1)})\,r(\psi^i_j\ox S(a_{(2)}))\,\gamma^j$, and
$\gamma^1,\,\ldots,\,\gamma^N$ is a basis of $\Gamma$ as a left
$\cB$-mod\-ule, too. \qed
\end{lemma}

We set $\cA=\cO(\SL_q(2))$ and denote by $u^i_j$, $i,\,j=1,\,2$, the
canonical generators of $\cA$. Then $r$ is defined by
$r(u^i_k\ox u^j_l) = q^{-1/2}\,R^{ij}_{kl}$, where $R$ is the R-ma\-trix of
$\SL_q(2)$ specified in \cite{RTF}.

If $c=-q/(q+1)^2$ and $x_i\in\Lin_\dC\{u^1_i,\,u^2_i\}$, $i=1,\,2$,
with
\begin{displaymath}
\frac{\e(x_1)}{\e(x_2)} = \frac{q+1}{q}\,\frac{\alpha_1}{\alpha_0}
\quad \mbox{and} \quad \e(x_1)^2
= -\frac{q+1}{q\,(q-1)}\,\frac{\alpha_1}{\alpha_0},
\end{displaymath}
then the subalgebra $\cB$ of $\cA$ generated by $x_1$, $x_2$ is
generated by $x_1$, $x_2$ with the relation
$x_1\,x_2 - q\,x_2\,x_1 = 1$, and it contains $\cB_c$ as the
subalgebra of the elements of even degree ($x_1$, $x_2$ being of
degree 1), cf.~\cite{Mue1}. The lemma with $N=2$, $b_i=x_i$ and
$\nu=\id$ yields an equivariant derivation $\de$ of $\cB$, and
$\de|_{\cB_c}$ is an equivariant derivation of $\cB_c$, which is by
construction independent of the embedding of $\cB_c$ in $\cA$. We
obtain
\begin{displaymath}
\textstyle \de a\,b = \sum_i\,\gamma^i\,
\bigl(x_i\,a - \sum_j\,a_{(1)}\,x_j\,r(u^j_i\ox a_{(2)})\bigr)\,b
\end{displaymath}
for all $a,\,b\in\cB$. In particular,
$\omega = \gamma^1\,x_1 + \gamma^2\,x_2$,
\begin{displaymath}
\begin{array}{ll}
\de e_{-1} & = (q-1)\,(-\omega\,e_{-1}
+ q^{-2}\,\alpha_0\,\gamma^1\,x_2), \\*[\xx]
\de e_0 & = (q-1)\,(-\omega\,e_0
- q^{-2}\,\alpha_0\,(\gamma^1\,x_1 - q^2\,\gamma^2\,x_2)), \\*[\xx]
\de e_1 & = (q-1)\,(-\omega\,e_1
+ \alpha_0\,\gamma^2\,x_1)
\end{array}
\end{displaymath}
and $\omega = \frac{q^2\,(q+1)^2}{(q-1)\,(q^3-1)\,\alpha_0^2}\,
\bigl(\de e_{-1}\,e_1 + (q^2+1)^{-1}\,\de e_0\,e_0
+ q^{-2}\,\de e_1\,e_{-1}\bigr)$, thus
$\gamma^i\,x_j\in\Gamma(\de|_{\cB_c})$ for all $i,\,j=1,\,2$. We set
$\gamma^{in} = \gamma^i\,(-q)^{\delta_{n1}}\,x_{3-n}$ and calculate that
\begin{displaymath}
\textstyle \de a\,b = \sum_{in}\,\gamma^{in}\,x_n\,
\bigl(x_i\,a - \sum_j\,a_{(1)}\,x_j\,r(u^j_i\ox a_{(2)})\bigr)\,b
\end{displaymath}
for all $a,\,b\in\cB$, thus $\overline{\de a\,b}
= \overline{\sum_{in}\,\gamma^{in}\,\e(x_n)\,\chi_i(a)\,\e(b)}$ for all
$a,\,b\in\cB_c$, where $\overline{\hspace{0.75em}\rule{0pt}{1.25ex}}:
\Gamma(\de|_{\cB_c})\to\Gamma(\de|_{\cB_c})/(\Gamma(\de|_{\cB_c})\cB_c^+)$
is the canonical projection and $\chi_i(a) = \e(x_i)\,\e(a) - r(x_i\ox a)$.
Since $\dim_\dC\Lin_\dC\{\chi_1|_{\cB_c},\,\chi_2|_{\cB_c}\} = 2$ and
$\dim_\dC\Lin_\dC\bigl\{\overline{\sum_n\,\gamma^{in}\,\e(x_n)}
\;\big|\;i=1,\,2\bigr\} = 2$, this implies
$\dim_\dC\overline{\Gamma(\de|_{\cB_c})} = 2$, that is,
$\dim_{\e,{\rm r}}(\de|_{\cB_c}) = 2$. Our proof also shows that
$\{\de e_{-1},\,\de e_0,\,\de e_1\}$ generates $\Gamma(\de|_{\cB_c})$
as a right $\cB_c$-mod\-ule. A similar argument starting with
\begin{displaymath}
\begin{array}{rl}
a\,\de b & = \sum_{ijkl}
\,a\,\bigl(x_l\,b_{(1)}\,r(u^k_i\ox S(b_{(2)})) - b\,x_l\,\delta_{ik}\bigr)
\,r(u^j_k\ox S(u^l_j))\,\gamma^i \\*[\xx]
& = \sum_{ijk}
\,a\,\bigl(\sum_l\,x_l\,b_{(1)}\,r(S(u^l_j)\ox b_{(2)}) - b\,x_j\bigr)
\,r(u^k_i\ox S(u^j_k))\,\gamma^i
\end{array}
\end{displaymath}
shows $\dim_{\e,{\rm l}}(\de|_{\cB_c}) = 2$. For the definition of
$\cO(\Sph^2_{qc})$ as a right $\cO(\SO_{q^2}(3))$-co\-mod\-ule algebra
only the square of $q$ is needed, therefore we can replace $q$ by $-q$
and obtain the corresponding result for $c=q/(-q+1)^2$. This
establishes the claims about the solutions 1 and 2 in
Section \ref{classif}.

If $c=c(2)$ and $x_i,\,y_i\in\Lin_\dC\{u^1_i,\,u^2_i\}$, $i=1,\,2$, with
either
\begin{displaymath}
\begin{array}{l}
\mbox{(i)} \quad \frac{\e(x_1)}{\e(x_2)}
= \frac{(q^4+1)\,\alpha_1}{q^4\,\alpha_0}
\quad \mbox{and} \quad \frac{\e(y_1)}{\e(y_2)}
= \frac{(q^4+1)\,\alpha_1}{q\,\alpha_0} \qquad \mbox{or} \\*[\xx]
\mbox{(ii)} \quad \frac{\e(x_1)}{\e(x_2)}
= \frac{(q^4+1)\,\alpha_1}{\alpha_0}
\quad \mbox{and} \quad \frac{\e(y_1)}{\e(y_2)}
= \frac{(q^4+1)\,\alpha_1}{q^5\,\alpha_0},
\end{array}
\end{displaymath}
then $y_i\,x_j\in\cB_c$ for all $i,\,j=1,\,2$, and
$x_1\,y_2 - q\,x_2\,y_1 = \zeta\,1$ with
$\zeta\in\dC\setminus\{0\}$. Moreover, $\sum_m\,u^i_m\,\tau_{jm}(e_k)
= \sum_{mn}\,\tau_{mi}(e_n)\,u^m_j\,\pi^n_k$ for the solution (f) in
Section \ref{classif}. We consider the $\cB_c$-$\cA$-bi\-mod\-ule with
the right $\cA$-mod\-ule basis $\gamma^1$, $\gamma^2$ and the left
$\cB_c$-mod\-ule operation
$a\,\gamma^j = \sum_i\,\tau_{ij}(a)\,\gamma^i$. Then
$\de a := \omega\,a - a\,\omega$ with $\omega = \sum_i\,\gamma^i\,x_i$
defines an equivariant derivation of $\cB_c$:
\begin{displaymath}
\textstyle \de a\,b = \sum_i\,\gamma^i\,
\bigl(x_i\,a - \sum_j\,\tau_{ij}(a)\,x_j\bigr)\,b
\end{displaymath}
for all $a,\,b\in\cB_c$. In particular,
$\omega = \gamma^1\,x_1 + \gamma^2\,x_2$,
\begin{displaymath}
\begin{array}{ll}
\de e_{-1} & = \omega\,e_{-1} - \frac{q^2-1}{q^4+1}\,\alpha_0\,
\gamma^1\,x_2, \\*[\xx]
\de e_0 & = \omega\,e_0 + \frac{q^2-1}{q^4+1}\,\alpha_0\,
(\gamma^1\,x_1 - q^2\,\gamma^2\,x_2), \\*[\xx]
\de e_1 & = \omega\,e_1 - q^2\,\frac{q^2-1}{q^4+1}\,\alpha_0\,
\gamma^2\,x_1.
\end{array}
\end{displaymath}
In addition,
$\omega = \frac{(q^4+1)^2}{(q^2-1)^2\,(q^2+1)\,\alpha_0^2}\,
\bigl(\de e_{-1}\,e_1 + (q^2+1)^{-1}\,\de e_0\,e_0
+ q^{-2}\,\de e_1\,e_{-1}\bigr)$ in the case (ii), in which we see
that $\{\de e_{-1},\,\de e_0,\,\de e_1\}$ generates $\Gamma(\de)$ as a
right $\cB_c$-mod\-ule. Using
$\gamma^{in} = \gamma^i\,\zeta^{-1}\,(-q)^{\delta_{n1}}\,x_{3-n}$, we
can proceed like before to show $\dim_{\e,{\rm r}}\de = 2$. In the
case (i), since the left action of
$\Lin_\dC\{1,\,e_{-1},\,e_0,\,e_1\}$ on
$\Lin_\dC\{\gamma^1,\,\gamma^2\}$ consists of all $\dC$-lin\-e\-ar
endomorphisms, $\{\de e_j,\,e_i\,\de e_j\;|\;i,\,j=-1,\,0,\,1\}$
generates $\Gamma(\de)$ as a right $\cB_c$-mod\-ule. Exploiting
e.g.~the relations
\begin{displaymath}
\begin{array}{c}
q^2\,e_{-1}\,\de e_0 =
\bigl(e_0 - q^2\,\frac{q^2-1}{q^4+1}\,\alpha_0\bigr)\,\de e_{-1},
\quad q^{-2}\,e_1\,\de e_0 =
\bigl(e_0 + \frac{q^2-1}{q^4+1}\,\alpha_0\bigr)\,\de e_1, \\*[\xx]
\de e_{-1}\,\bigl(q^2\,e_0 - \frac{q^6-1}{q^4+1}\,\alpha_0\bigr)
= \de e_0\,e_{-1},
\quad \de e_1\,\bigl(e_0 + \frac{q^6-1}{q^4+1}\,\alpha_0\bigr)
= q^2\,\de e_0\,e_1, \\*[\xx]
(q^2+1)\,\de e_{-1}\,e_1 + \de e_0\,e_0
+ (q^{-2}+1)\,\de e_1\,e_{-1} = 0
\end{array}
\end{displaymath}
we get $\dim_{\e,{\rm r}}\de \leq 2$, while
$\dim_{\e,{\rm r}}\de \geq 2$ follows from Theorem \ref{th1} (iii). 
However, one checks that in both cases $\dim_{\e,{\rm l}}\de = 0$. We
identify the cases (i) and (ii) with the solutions 3 and 4 in Section
\ref{classif}.

We take the basis $\tilde e_{-2},\,\ldots,\,\tilde e_2$ in \cite{Po2}
of the spin 2 subcomodule of $\cB_c$ for $b_1,\,\ldots,\,b_N$ and set
$\nu=\id$ in the above lemma. This yields the equivariant derivation
$\de$ of $\cB_c$ which is given by
$\de b\,a = \tilde\chi(b_{(1)})\ox b_{(2)}\,a$ with
$\tilde\chi_i(b) = r(\tilde e_i\ox b) - \e(\tilde e_i)\,\e(b)$ for all
$a,\,b\in\cB_c$, $i=-2,\,\ldots,\,2$. If $c=0$, then it corresponds to
each case of the solution 5 in Section \ref{classif}. To compute this,
we write $\tilde\chi_i$ in terms of $l(a) = r(u^1_1\ox a)$,
$l^{-1}(a) = r(u^2_2\ox a)$ and $l^{(-)}(a) = -q\,r(u^2_1\ox a)$,
while $r(u^1_2\ox a) = 0$, and use the relations
$l\,l^{-1} = \e = l^{-1}\,l$, $l\,l^{(-)} = q\,l^{(-)}\,l$ and
$(f\,g)|_{\cB_c} = 0$, where $f$ is defined in Section \ref{classif}
and $g$ is any $\dC$-lin\-e\-ar functional on $\cA$. Hence, we can
restrict ourselves to $\alpha_{-1}=0=\alpha_1$ to calculate
\begin{displaymath}
\de(e_{-1}^3) = q^{-6}\,(q^4+q^2+1)\,
(q^2\,\de\tilde e_{-2}\,e_{-1} - \de e_{-1}\,\tilde e_{-2}).
\end{displaymath}
Since $\de$ is equivariant, this implies that for each $a$ in the
spin 3 subcomodule of $\cB_c$, $\de a$, and thus $\Gamma(\de)$, is
contained in the right $\cB_c$-mod\-ule generated by
$\{\de e_{-1},\,\de e_0,\,\de e_1,\,
\de\tilde e_{-2},\,\ldots,\,\de\tilde e_2\}$. Exploiting e.g.~the
relations
\begin{displaymath}
\begin{array}{c}
(q^2+1)\,\de e_{-1}\,e_1 + \de e_0\,e_0
+ (q^{-2}+1)\,\de e_1\,e_{-1} = 0, \\*[\xx]
\de e_{-1}\,(q^2\,e_0 + \alpha_0) = \de e_0\,e_{-1}, \\*[\xx]
(q^4+1)\,\de\tilde e_{-2}\,(q^2\,e_0 + \alpha_0)
- \de\tilde e_{-1}\,e_{-1}
= q^4\,(q^2+1)\,\alpha_0\,\de e_{-1}\,e_{-1}, \\*[\xx]
\alpha_0^2\,\de e_{-1} = \frac{q^4+1}{q^2}\,\de\tilde e_{-2}\,e_1
+ \frac{1}{q^2\,(q^2+1)}\,\de\tilde e_{-1}\,e_0
- \frac{1}{q^2\,(q^4+q^2+1)}\,\de\tilde e_0\,e_{-1}, \\*[\xx]
\alpha_0^2\,\de e_0 = \de\tilde e_{-1}\,e_1
+ \frac{q^4-1}{q^6-1}\,\de\tilde e_0\,e_0
+ q^{-4}\,\de\tilde e_1\,e_{-1}, \\*[\xx]
\alpha_0^2\,\de e_1 = -\frac{q^4}{q^4+q^2+1}\,\de\tilde e_0\,e_1
+ \frac{1}{q^2+1}\,\de\tilde e_1\,e_0
+ \frac{q^4+1}{q^4}\,\de\tilde e_2\,e_{-1}
\end{array}
\end{displaymath}
we get $\dim_{\e,{\rm r}}\de \leq 2$, while
$\dim_{\e,{\rm r}}\de \geq 2$ follows from Theorem \ref{th1} (iii).
Similarly, the derivation associated to the solution 6, say, $\de'$,
arises if $r'$ with $r'(a\ox b) = r(S(b)\ox a)$, $a,\,b\in\cA$, is
used instead of $r$. With regard to the left $\cB_c$-mod\-ule
structure of $\Gamma$ the derivation $\de$ is given by
$a\,\de b = \tilde\chi'(b_{(1)})\ox a\,b_{(2)}$ with
$\tilde\chi_i'(b) = r(S(\tilde e_i)\ox b) - \e(\tilde e_i)\,\e(b)$ for
all $a,\,b\in\cB_c$, $i=-2,\,\ldots,\,2$, thus
$\de^*b\,a = \tilde\chi'(b_{(1)}^*)\ox (b_{(2)}\,a)^*$.
One checks that $\de'=\de^*$, if $e_i^*=e_{-i}$, and therefore
$\dim_{\e,{\rm l}}\de = \dim_{\e,{\rm r}}\de^* = 2$. Finally, the
2-di\-men\-sional covariant differential calculus described in
\cite{Po3} corresponds to the solution~7. It is equal to
$\de_{(\cB_c^+)^2}$, if $\alpha_{-1}=0=\alpha_1$.

The classification problem in Section \ref{classif} with
$\dim_{\e,{\rm l}}\de = 2$ instead of $\dim_{\e,{\rm r}}\de = 2$ is
equivalent, because the Hopf algebra $\cO(\SL_q(2))$ with the opposite
multiplication is isomorphic to $\cO(\SL_{q^{-1}}(2))$ and
correspondingly the right comodule algebra $\cO(\Sph^2_{qc})$ with the
opposite multiplication to $\cO(\Sph^2_{q^{-1} c})$. It is still open
for $c=c(0),\,c(1),\,\ldots\,$ and for $\rho=\lambda=0$. Since in the
one-to-one correspondence in Theorem \ref{th1} (iii) only the trivial
derivation $\de=\de_{\cB^+}$ satisfies $\dim_{\e,{\rm l}}\de = 0$, the
equivalents with the opposite comultiplication of the equivariant
derivations for $c=c(2)$ constructed above do not occur, i.e.~they do
not arise from right ideals.

\section*{\rm ACKNOWLEDGMENTS}

I am grateful to Prof.~K.~Schm\"udgen for showing his interest in my
work and to S.~Kolb for detailed discussions. This work was supported
by the Deut\-sche For\-schungs\-ge\-mein\-schaft with\-in the scope of
the postgraduate scholarship programme ``Gra\-du\-ier\-ten\-kol\-leg
Quan\-ten\-feld\-theo\-rie'' at the University of Leipzig.

\def\refname


{\rm REFERENCES}
\begin{thebibliography}{99}

\bibitem{AS} J.~Apel and K.~Schm\"udgen, Classification of
three-di\-men\-sional covariant differential calculi on Podles' quantum
spheres and on related spaces, {\em Lett.\ Math.\ Phys.\/} {\bf 32} (1994),
25--36.

\bibitem{Di} M.~S.\ Dijkhuizen, Some remarks on the construction of quantum
symmetric spaces, {\em Acta Appl.\ Math.\/} {\bf 44} (1996), 59--80,
\mbox{\tt math.QA/9512225 \hspace{-0.5em}}.

\bibitem{DK} M.~S.\ Dijkhuizen and T.~H.\ Koornwinder, Quantum homogeneous
spaces, duality, and quantum 2-spheres, {\em Geom.\ Dedicata\/} {\bf 52}
(1994), 291--315.

\bibitem{Her} U.~Hermisson, Construction of covariant differential calculi
on quantum homogeneous spaces, {\em Lett.\ Math.\ Phys.\/} {\bf 46} (1998),
313--322, \mbox{\tt math.QA/9806008 \hspace{-0.5em}}.

\bibitem{Mue1} E.~F.\ M\"uller, Konstruktion von Rechtscoidealunteralgebren,
Degree Dissertation, Munich, 1995.

\bibitem{Mue2} E.~F.\ M\"uller, private communication, October 1998.

\bibitem{MS} E.~F.\ M\"uller and H.-J.\ Schneider, Quantum homogeneous
spaces with faithfully flat module structures, {\em Israel J.\ Math.\/}
{\bf 111} (1999), 157--190.

\bibitem{NS} M.~Noumi and T.~Sugitani, Quantum symmetric spaces and related
q-orthogonal polynomials, in: A.\ Arima et al.\ (eds), {\em Group
theoretical methods in Physics,\/} Proceedings XX ICGTMP, Toyonaka (Japan),
1994, World Scientific, Singapore, 1995, p.~28--40,
\mbox{\tt math/9503225 \hspace{-0.5em}}. 

\bibitem{Po1} P.~Podle\'s, Quantum spheres, {\em Lett.\ Math.\ Phys.\/}
{\bf 14} (1987), 193--202.

\bibitem{Po2} P.~Podle\'s, Differential calculus on quantum spheres,
{\em Lett.\ Math.\ Phys.\/} {\bf 18} (1989), 107--119.

\bibitem{Po3} P.~Podle\'s, The classification of differential structures on
quantum 2-spheres, {\em Comm.\ Math.\ Phys.\/} {\bf 150} (1992), 167--179.

\bibitem{RTF} N.~Yu.\ Reshetikhin, L.~A.\ Takhtadzhyan and L.~D.\ Faddeev,
Quantization of Lie groups and Lie algebras, {\em Leningrad Math.\ J.\/}
{\bf 1} (1990), 193--225.

\bibitem{SD} J.~V.\ Stokman and M.~S.\ Dijkhuizen, Quantized flag manifolds
and irreducible $*$-rep\-re\-sen\-ta\-tions, {\em Comm.\ Math.\ Phys.\/}
{\bf 203} (1999), 297--324, \mbox{\tt math.QA/9802086 \hspace{-0.5em}}.

\bibitem{Sw} M.~E.\ Sweedler, {\em Hopf algebras,\/} W.~A.\ Benjamin, New
York, 1969.

\bibitem{T} M.~Takeuchi, Relative Hopf modules---equivalences and freeness
criteria, {\em J.\ Algebra\/} {\bf 60} (1979), 452--471.

\bibitem{Wo} S.~L.\ Woronowicz, Differential calculus on compact matrix
pseudogroups (quantum groups), {\em Comm.\ Math.\ Phys.\/} {\bf 122} (1989),
125--170.

\end{thebibliography}
\end{document}